\documentclass[3pt]{article}
\usepackage{graphics, amssymb, amsmath, mathrsfs,graphicx,theorem}
\usepackage[numbers,sort&compress]{natbib}
\newtheorem{theorem}{Theorem}[section]

\newtheorem{lemma}{Lemma}[section]

\newtheorem{definition}{Definition}[section]
\newtheorem{remark}{Remark}[section]

\numberwithin{equation}{section}

\begin{document}

\title{A new study on the mild solution for impulsive fractional evolution equations}

\author{ Xiao-Bao Shu$^{a}$
\thanks{ Corresponding author: Xiao-Bao Shu.
Email: sxb0221@163.com(Xiao-Bao Shu),
shulinxin2066@126.com(Linxin Shu),
fxu.feixu@gmail.com(Fei Xu).
}
, Linxin Shu$^{a}$
, Fei Xu$^{b}$
\\
\footnotesize{
$^a$ College of Mathematics and Econometrics, Hunan University,}
\footnotesize{Changsha, Hunan 410082, PR China.}\\
\footnotesize{$^b$ Department of Mathematics, Wilfrid Laurier University,}
\footnotesize{Waterloo, Ontario, N2L 3C5, Canada.}\\
}

\date{}
\maketitle
\begin{center}
\begin{minipage}{139mm}

\noindent{\bf Abstract}
In this article, we consider mild solutions to a class of impulsive fractional evolution equations of order $0<\alpha<1$. After analyzing analytic results reported in the literature using Mittag-Leffer function, $\alpha$-resolvent operator theory, we propose a more appropriate new definition of mild solutions for impulsive fractional evolution equations by replacing  the impulse term operator  $S_\alpha(t-t_i)$ with $S_\alpha(t)S_\alpha^{-1}(t_i)$, where $S_\alpha^{-1}(t_i)$ denotes the inverse of the fractional solution operator $S_\alpha(t)$ at $t=t_i,~(i=1,2,\cdots m)$.

\end{minipage}
\end{center}\vspace{5mm}

\section{Introduction \label {Sec1}}
Fractional differential equations have been used to describe the dynamics of a variety of engineering systems and financial systems \cite{5,shu12,E function,s operator,up,GYC,XF}. Impulsive systems are an active research area and have drawn more and more attention of researchers in mathematics and applied fields \cite{n1,n2,n3,6,7,8,9,10,DSF,S.L2018}.
In particular, impulsive fractional evolution equations recently received considerable attention in the literature.
The existence, uniqueness and other properties of the mild solutions to impulsive fractional evolution equations have been investigated in many works \cite{n1,n2,n3,6,7,8,9,10,11,12,13,14,15}.
Recently, the form of solutions to impulsive fractional evolution equations was studied \cite{2,4,1}.

As an example, we consider the following general impulsive fractional evolution equation with order $0<\alpha<1$.
\begin{equation}\label{eq1}
\left\{ \begin{array}{ll}
_0D_t^\alpha x(t)=Ax(t)+f(t,x(t)),\;\; 0<\alpha<1,\;t\in J=[0,T],\;t\ne t_k,\\
\Delta x|_{t=t_k}=I_k(x(t_k^-)),\;(k=1,2,\ldots,m,),\\
x(0)=x_0,\\
\end{array} \right.
\end{equation}
where $_0D_t^\alpha$ is Caputo fractional derivative, $A$ is a sectorial operator of type $(M,\theta,\alpha,\mu)$ in complex Banach space $X$,
$f:J\times X\rightarrow X$ is a continuous function,
$I_k:X\rightarrow X (k=1,2,\cdots ,m)$ are appropriate functions, and $0=t_0<t_1<t_2<\cdots<t_m<t_{m+1}=T$. Let $J_k=(t_k,t_{k+1}], k=1,2,\cdots,m,$ $J_0=[0,t_1],$ and  $\triangle x|_{t = t_k }=x(t_k^+)-x(t_k^-)$, where $x(t_k^+)$ and $x(t_k^-)$ denote the right and the left limits of $x(t)$ at $t=t_k$, respectively.

The mild solutions to  equation \eqref{eq1} were investigated in the literature and two types of solutions are obtained.

The first type of  solution is obtained by taking integrals over $(t_k,\,t_{k+1}]$, $(k=1,2,  \cdots ,m)$ and $[0,\,t_1]$ , given by \cite{2,n1,n2,n3}
\begin{equation}\label{sol1}
      x(t)=\left\{ \begin{array}{ll}
         S_\alpha(t)x_0+\int_{0}^tT_\alpha(t-\theta)f(\theta,x(\theta))d\theta,  & t\in[0,t_1],\\
         S_\alpha(t-t_1)[x(t_1^-)+I_1(x(t_1^-))]
          +\int_{t_1}^tT_\alpha(t-\theta)f(\theta,x(\theta))d\theta,\;&t\in(t_1,t_{2}],\\
         \vdots\\
         S_\alpha(t-t_m)[x(t_m^-)+I_m(x(t_m^-))]
          +\int_{t_m}^tT_\alpha(t-\theta)f(\theta,x(\theta))d\theta,\;&t\in(t_m,T].
        \end{array} \right.
 \end{equation}

The second type of solution is obtained by taking integrals over $[0,t]$, given by \cite{4,1,11,12,13,14,15}
\begin{equation}\label{sol2}
      x(t)=\left\{ \begin{array}{ll}
         S_\alpha(t)x_0+\int_{0}^tT_\alpha(t-\theta)f(\theta,x(\theta))d\theta,  & t\in[0,t_1],\\
          S_\alpha(t)x_0+\sum\limits_{i=1}^k S_\alpha(t-t_i)I_i(x(t_i^-))
          +\int_{0}^tT_\alpha(t-\theta)f(\theta,x(\theta))d\theta,\;&t\in(t_k,t_{k+1}].
        \end{array} \right.
 \end{equation}

Recently, Shu et al.\cite{1} proved that it is inappropriate to use integrals over $[0,\,t_1]$ and $(t_k,\,t_{k+1}]$, $(k=1,2,  \cdots ,m)$ to describe mild solutions to impulsive fractional evolution equation when the fractional derivative is $_{0}D^\alpha_t$ since
\begin{equation}\label{DtT}
_{t_k}D^\alpha_t\Big(\int_{t_k}^tT_\alpha(t-\theta)f(\theta)d\theta\Big)
  =A\int_{t_k}^t T_\alpha(t-\theta)f(\theta)d\theta+f(t),
\end{equation}
and
\begin{equation}\label{D0Dt}
_{0}D^\alpha_t\Big(\int_{t_k}^tT_\alpha(t-\theta)f(\theta)d\theta\Big)
\ne~_{0}D^\alpha_t\Big(\int_{0}^tT_\alpha(t-\theta)f(\theta)d\theta\Big).
\end{equation}
It thus follows that the integral interval of mild solutions should be consistent with the integral interval of fractional derivative operator. We find that there exist similar problems in (1.2) and (1.3) , when impulsive term operator $S_\alpha(t-t_i)I_i(x(t_i^-)),(i=1,2,\cdots,m)$ under the fractional derivative operate $_0D_t^\alpha$, i.e. $$_{t_i}D^\alpha_t\Big(S_\alpha(t-t_i)I_i(x(t_i^-))\Big)= AS_\alpha(t-t_i)I_i(x(t_i^-)),~~
_{0}D^\alpha_t\Big(S_\alpha(t-t_i)I_i(x(t_i^-))\Big)\neq AS_\alpha(t-t_i)I_i(x(t_i^-)).$$
We shall show these in Lemma \ref{lem2} and Remark \ref{remA}.
Thus, the expressions of mild solution (1.2) and (1.3) are both  incorrect. To the best of our knowledge, appropriate expressions of mild solution to impulsive  fractional evolution equations of order $0<\alpha<1$ have not been established in the literature.

In this paper, we study system \eqref{eq1} and give a correct solution to the system as shown in Theorem 1.1.  Here, instead of using $\sum\limits_{i=1}^k S_\alpha(t-t_i)I_i(x(t_i^-)))$, we use $\sum\limits_{i=1}^k S_\alpha(t)S^{-1}_\alpha(t_i)I_i(x(t_i^-)))$ for the impulse terms of the solution to avoid the above mentioned problem.

\begin{theorem}\label{thm3}
A function $x(t)\in PC(J,X)$ is called a mild solution of problem \eqref{eq1}, if it satisfies the integral equation
\begin{equation}\label{sol3}
      x(t)=\left\{\begin{array}{ll}
         S_\alpha(t)x_0+\int_{0}^tT_\alpha(t-\theta)f(\theta,x(\theta))d\theta,  & t\in[0,t_1],\\
          S_\alpha(t)x_0+ S_\alpha(t)\sum\limits_{i=1}^k S^{-1}_\alpha(t_i)I_i(x(t_i^-))
          +\int_{0}^tT_\alpha(t-\theta)f(\theta,x(\theta))d\theta,\;&t\in(t_k,t_{k+1}].
        \end{array} \right.
 \end{equation}
\end{theorem}

The rest of this paper is organized as follows: In Section 2, we introduce some fundamental definitions and lemmas. In Section 3, we give a linear impulsive fractional differential equation as an example to illustrate that the above two kinds of solution forms are not correct.
The proofs of Theorems 1.1 will be given in Section 4.

\section{Preliminaries \label {Sec2}}
In this section, we will give some necessary definitions and lemmas.

\begin{definition}
The Caputo derivative of order $\alpha$ for a function  $f : [a, +\infty)  \rightarrow  X$ is
defined as
$$ _aD^{\alpha}_t f(t)=\frac{1}{\Gamma(m-\alpha)}\int_{a}^{t}(t-\tau)^{m-\alpha-1}f^{(m)}(\tau)d\tau, \;\;\;m-1<\alpha\leq m, t\in [a,+\infty).$$
\end{definition}

\begin{definition}\cite{E function}
The Mittag-Leffer function is defined as
$$ E_{{\alpha,\beta}}   \left( z \right) =\sum _{j=0}^{
\infty }{\frac {{z}^{j}}{\Gamma  \left( \alpha\,j+\beta \right) }}={\frac {1}{2 \, \pi \,i}} \int_{H_{\alpha}} {\frac {{e}^{\mu}{\mu}^{\alpha-\beta}}{{\mu}^{\alpha}-z}} d \mu,\quad \alpha,\,\beta\,>0,\,z \in \mathbb{C},
$$
where $H_{\alpha}$ is a Hankel path, a contour starting and ending at $-\infty$, and
encircling the disc $|\mu|<|z|^{\frac{1}{\alpha}}$  counterclockwise.
\end{definition}

\begin{definition}\cite{2}
Let $A: D\subseteq X \to X$ be a closed linear operator.
$A$ is said to be sectorial operator of type $(M,\theta, \alpha, \mu) $.
If there exist $0<\theta<\pi/2$, $M>0$
and $\mu\in \mathbb{R}$ such that the $\alpha$-resolvent of $A$ exists
outside the sector
\[
\mu+ S_\theta=\{\mu+\lambda^\alpha: \lambda\in \mathbb{C},
|\arg(-\lambda^\alpha)|<\theta\} \label {eqn2.1}
\]
and
\[
\|(\lambda^\alpha I-A)^{-1}\|\leq\frac{M}{|\lambda^\alpha-\mu|}, \quad \lambda^\alpha
\not \in \mu+ S_\theta.  \label {eqn2.2}
\]
\end{definition}

\begin{lemma}\label{lem0}\cite{1}\cite{up}
If $A$ is a sectorial operator of type $(M,\theta,\alpha,\mu),$ then it is not difficult to see that $A$ is the infinitesimal generator of the $\alpha$-resolvent families $ \{S_\alpha(t)\}_{t\geq 0}, \{T_\alpha(t)\}_{t\geq 0}$,
\begin{align}\label{STK}
S_\alpha(t)&=\frac{1}{2\pi i}\int_c e^{\lambda t}\lambda^{\alpha-1}R(\lambda^\alpha,A)d\lambda=E_{\alpha,1}(At^\alpha)
=\sum_{k=0}^\infty\frac{(At^\alpha)^k}{\Gamma(1+\alpha k)},\\
T_\alpha(t)&=\frac{1}{2\pi i}\int_c e^{\lambda t}R(\lambda^\alpha,A)d\lambda=t^{\alpha-1}E_{\alpha,\alpha}(At^\alpha)
=t^{\alpha-1}\sum_{k=0}^\infty\frac{(At^\alpha)^k}{\Gamma(\alpha+\alpha k)},
\end{align}
where $c$ is a suitable path.
\end{lemma}

\begin{lemma}\label{lem1}\cite{1}
Suppose $A$ is a sectorial operator of type $(M,\theta,\alpha,\mu)$.
If  $0<\alpha<1$, then
\begin{equation}\label{lem1-1}
_0D^\alpha_t[S_\alpha(t)x_0]=A[S_\alpha(t)x_0],
\end{equation}
and
\begin{equation}\label{lem1-2}
 _0D^\alpha_t\big(\int_0^tT_\alpha(t-\theta)f(\theta)d\theta\big)
=A\int_0^t T_\alpha(t-\theta)f(\theta)d\theta+f(t).
\end{equation}
\end{lemma}

\section{A counter example \label {Sec3}}
In this section, we present a counter example to show that  \eqref{sol1} and \eqref{sol2} are incorrect. Consider the following linear Caputo fractional differential equation:
\begin{equation}\label{eq2}
\left\{ \begin{array}{ll}
_0D_t^{\frac{2}{3}} x(t)=\rho x(t)+t,\;\;\;t\in J=[0,T],\;t\ne t_1,\\
\Delta x|_{t=t_1}=y_1\in X,\\
x(0)=x_0.\\
\end{array} \right.
\end{equation}

In the following, considering  \eqref{sol1} and \eqref{sol2}, we give two types of  classical solutions to system \eqref{eq2}, respectively,  using Mittag-Lefller functions
\begin{equation}\label{sol01}
      x(t)=\left\{ \begin{array}{ll}
      E_{\frac{2}{3},1}(\rho t^{\frac{2}{3}})x_0 +\int_{0}^t(t-\theta)^{-\frac{1}{3}}E_{\frac{2}{3},\frac{2}{3}}(\rho(t-\theta)^{\frac{2}{3}})\theta d\theta,  &t\in[0,t_1],\\
      E_{\frac{2}{3},1}(\rho (t-t_1)^{\frac{2}{3}})[x(t_1^-)+y_1]
      +\int_{t_1}^t(t-\theta)^{-\frac{1}{3}}E_{\frac{2}{3},\frac{2}{3}}(\rho(t-\theta)^{\frac{2}{3}})\theta d\theta,
      &t\in(t_1,T],
      \end{array} \right.
 \end{equation}
 and\vspace*{-8pt}
\begin{equation}\label{sol02}
      x(t)=\left\{ \begin{array}{ll}
      E_{\frac{2}{3},1}(\rho t^{\frac{2}{3}})x_0 +\int_{0}^t(t-\theta)^{-\frac{1}{3}}E_{\frac{2}{3},\frac{2}{3}}(\rho(t-\theta)^{\frac{2}{3}})\theta d\theta,  &t\in[0,t_1],\\
      E_{\frac{2}{3},1}(\rho t^{\frac{2}{3}})x_0
      +E_{\frac{2}{3},1}(\rho (t-t_1)^{\frac{2}{3}})y_1
      +\int_{0}^t(t-\theta)^{-\frac{1}{3}}E_{\frac{2}{3},\frac{2}{3}}(\rho(t-\theta)^{\frac{2}{3}})\theta d\theta,
      &t\in(t_1,T].\\
      \end{array} \right.
 \end{equation}
Now, we show that solutions \eqref{sol01} and \eqref{sol02} do not satisfy system \eqref{eq2}.

\begin{theorem}\label{thm1}
The solution \eqref{sol01} does not satisfy system \eqref{eq2}.
\end{theorem}

{\bf Proof:}
It follows from \cite{s operator} that
\begin{align}
&_0D_t^\alpha E_{\alpha,1}(\rho t^\alpha)=\rho E_{\alpha,1}(\rho t^\alpha), ~~\alpha>0,\\
&_0D_t^\alpha \int_0^t(t-\theta)^{\alpha-1}E_{\alpha,\alpha}(\rho(t-\theta)^\alpha)f(\theta)d\theta
=\rho \int_0^t(t-\theta)^{\alpha-1}E_{\alpha,\alpha}(\rho(t-\theta)^\alpha)f(\theta)d\theta+f(t).
\end{align}

(i) For $t\in [0,t_1]$,
\vspace*{-5pt}
$$x(t)=E_{\frac{2}{3},1}(\rho t^{\frac{2}{3}})x_0 +\int_{0}^t(t-\theta)^{-\frac{1}{3}}E_{\frac{2}{3},\frac{2}{3}}(\rho(t-\theta)^{\frac{2}{3}})
\theta d\theta,~t\in [0,t_1],$$
combined with (3.4), (3.5), we have
\begin{align*}
 _0D_t^{\frac{2}{3}} x(t)&=~_0D_t^{\frac{2}{3}}\Big[E_{\frac{2}{3},1}(\rho t^{\frac{2}{3}})x_0 +\int_{0}^t(t-\theta)^{\frac{2}{3}-1}E_{\frac{2}{3},\frac{2}{3}}(\rho(t-\theta)^{\frac{2}{3}})\theta d\theta\Big] \\
 &=\rho\Big[E_{\frac{2}{3},1}(\rho t^{\frac{2}{3}})x_0 +\int_{0}^t(t-\theta)^{-\frac{1}{3}}E_{\frac{2}{3},\frac{2}{3}}(\rho(t-\theta)^{\frac{2}{3}})\theta d\theta\Big]+t
 =\rho x(t)+t,\;\;\;\;t\in [0,t_1],
\end{align*}
which satisfied the equation \eqref{eq2}.

(ii) For $t\in(t_1,T]$, by \eqref{sol01},
\vspace*{-5pt}
$$x(t)=E_{\frac{2}{3},1}(\rho (t-t_1)^{\frac{2}{3}})[x(t_1^-)+y_1]
      +\int_{t_1}^t(t-\theta)^{-\frac{1}{3}}E_{\frac{2}{3},\frac{2}{3}}(\rho(t-\theta)^{\frac{2}{3}})\theta d\theta,
      \;\;\;t\in(t_1,T].$$
From the definition of the Caputo derivative, together with (3.4), (3.5), we have
\begin{align}\label{C2.2}
_0D_t^{\frac{2}{3}} x(t)
&=\int_0^t\frac{(t-s)^{-\frac{2}{3}}}{\Gamma(\frac{1}{3})}E'_{\frac{2}{3},1}(\rho (s-t_1)^{\frac{2}{3}})[x(t_1^-)+y_1]ds\nonumber\\
&+\int_0^t\frac{(t-s)^{-\frac{2}{3}}}{\Gamma(\frac{1}{3})}
\Big[\int_{t_1}^s(s-\theta)^{-\frac{1}{3}}E_{\frac{2}{3},\frac{2}{3}}
(\rho(s-\theta)^{\frac{2}{3}})\theta d\theta \Big]'ds,~t\in(t_1,T],
\end{align}
where
\begin{align}\label{I1}
 &\int_0^t\frac{(t-s)^{-\frac{2}{3}}}{\Gamma(\frac{1}{3})}E'_{\frac{2}{3},1}(\rho (s-t_1)^{\frac{2}{3}})[x(t_1^-)+y_1]ds\nonumber\\
 &=[x(t_1^-)+y_1]\Big[\int_0^{t_1}\frac{(t-s)^{-\frac{2}{3}}}{\Gamma(\frac{1}{3})}E'_{\frac{2}{3},1}(\rho (s-t_1)^{\frac{2}{3}})ds
 +\int_{t_1}^{t}\frac{(t-s)^{-\frac{2}{3}}}{\Gamma(\frac{1}{3})}E'_{\frac{2}{3},1}(\rho (s-t_1)^{\frac{2}{3}})ds\Big]\nonumber\\
 &=[x(t_1^-)+y_1]\Big[\int_0^{t_1}\frac{(t-s)^{-\frac{2}{3}}}{\Gamma(\frac{1}{3})}E'_{\frac{2}{3},1}(\rho (s-t_1)^{\frac{2}{3}})ds
 +\int_{0}^{t-t_1}\frac{(t-t_1-s)^{-\frac{2}{3}}}{\Gamma(\frac{1}{3})}E'_{\frac{2}{3},1}(\rho s^{\frac{2}{3}})ds\Big]\nonumber\\
 &=[x(t_1^-)+y_1]\Big[\int_0^{t_1}\frac{(t-s)^{-\frac{2}{3}}}{\Gamma(\frac{1}{3})}E'_{\frac{2}{3},1}(\rho (s-t_1)^{\frac{2}{3}})ds+
 ~_0D_{t-t_1}^{\frac{2}{3}} E_{\frac{2}{3},1}(\rho t^{\frac{2}{3}})\Big]\nonumber\\
 &=[x(t_1^-)+y_1]\Big[\int_0^{t_1}\frac{(t-s)^{-\frac{2}{3}}}{\Gamma(\frac{1}{3})}E'_{\frac{2}{3},1}(\rho (s-t_1)^{\frac{2}{3}})ds+
 ~\rho E_{\frac{2}{3},1}(\rho (t-t_1)^{\frac{2}{3}})\Big],
\end{align}
and\vspace*{-5pt}
\begin{align}\label{I2}
&\int_0^t\frac{(t-s)^{-\frac{2}{3}}}{\Gamma(\frac{1}{3})}\Big[
\int_{t_1}^s(s-\theta)^{-\frac{1}{3}}E_{\frac{2}{3},\frac{2}{3}}
(\rho(s-\theta)^{\frac{2}{3}})\theta d\theta\Big]'ds\nonumber\\
&=\int_0^{t_1}\frac{(t-s)^{-\frac{2}{3}}}{\Gamma(\frac{1}{3})}\Big[
\int_{t_1}^s(s-\theta)^{-\frac{1}{3}}E_{\frac{2}{3},\frac{2}{3}}
(\rho(s-\theta)^{\frac{2}{3}})\theta d\theta\Big]' ds\nonumber\\
&+\int^t_{t_1}\frac{(t-s)^{-\frac{2}{3}}}{\Gamma(\frac{1}{3})}
\Big[\int_{t_1}^{s}(s-\theta)^{-\frac{1}{3}}E_{\frac{2}{3},\frac{2}{3}}
(\rho(s-\theta)^{\frac{2}{3}})\theta d\theta\Big]'ds\nonumber\\
&=\int_0^{t_1}\frac{(t-s)^{-\frac{2}{3}}}{\Gamma(\frac{1}{3})}\Big[
\int_{t_1}^s(s-\theta)^{-\frac{1}{3}}E_{\frac{2}{3},\frac{2}{3}}
(\rho(s-\theta)^{\frac{2}{3}})\theta d\theta\Big]'ds\nonumber\\
&+~_{t_1}D_t^{\frac{2}{3}}\int_{t_1}^{t}(t-\theta)^{-\frac{1}{3}}E_{\frac{2}{3},\frac{2}{3}}
(\rho(t-\theta)^{\frac{2}{3}})\theta d\theta\nonumber\\
&=\int_0^{t_1}\frac{(t-s)^{-\frac{2}{3}}}{\Gamma(\frac{1}{3})}\Big[
\int_{t_1}^s(s-\theta)^{-\frac{1}{3}}E_{\frac{2}{3},\frac{2}{3}}
(\rho(s-\theta)^{\frac{2}{3}})\theta d\theta\Big]'ds\nonumber\\
&+~\rho\int_{t_1}^{t}(t-\theta)^{-\frac{1}{3}}E_{\frac{2}{3},\frac{2}{3}}
(\rho(t-\theta)^{\frac{2}{3}})\theta d\theta+t.
\end{align}

Substituting \eqref{I1} and \eqref{I2} into \eqref{C2.2} yields
\begin{align}\label{1F(t)}
_0D_t^{\frac{2}{3}} x(t)
&=\rho\Big[E_{\frac{2}{3},1}(\rho(t-t_1)^{\frac{2}{3}})[x(t_1^-)+y_1]
    +\int_{t_1}^t(t-\theta)^{-\frac{1}{3}}E_{\frac{2}{3},\frac{2}{3}}
    (\rho(t-\theta)^{\frac{2}{3}})\theta d\theta\Big]+t\nonumber\\
&+[x(t_1^-)+y_1]\int_0^{t_1}\frac{(t-s)^{-\frac{2}{3}}}{\Gamma(\frac{1}{3})}E'_{\frac{2}{3},1}(\rho (s-t_1)^{\frac{2}{3}})ds\nonumber\\
&+\int_0^{t_1}\frac{(t-s)^{-\frac{2}{3}}}{\Gamma(\frac{1}{3})}\Big[
\int_{t_1}^s(s-\theta)^{-\frac{1}{3}}E_{\frac{2}{3},\frac{2}{3}}
(\rho(s-\theta)^{\frac{2}{3}})\theta d\theta\Big]' ds, t\in(t_1,T].
\end{align}
Let
\begin{align*}
&[x(t_1^-)+y_1]\int_0^{t_1}\frac{(t-s)^{-\frac{2}{3}}}{\Gamma(\frac{1}{3})}E'_{\frac{2}{3},1}(\rho (s-t_1)^{\frac{2}{3}})ds\nonumber\\
&+\int_0^{t_1}\frac{(t-s)^{-\frac{2}{3}}}{\Gamma(\frac{1}{3})}\Big[
\int_{t_1}^s(s-\theta)^{-\frac{1}{3}}E_{\frac{2}{3},\frac{2}{3}}
(\rho(s-\theta)^{\frac{2}{3}})\theta d\theta\Big]' ds=F(t),
\end{align*}
from this,
\begin{align*}
_0D_t^{\frac{2}{3}} x(t)=\rho x(t)+t+F(t),~ t\in(t_1,T].
\end{align*}

If we want the solution $x(t)$ given by \eqref{sol01} for $t\in (t_1,T]$ satisfies system \eqref{eq2}, i.e. $_0D_t^{\frac{2}{3}} x(t)=\rho x(t)+t$,
then $F(t)\equiv 0$ is valid for arbitrary $t\in (t_1,T]$,
However, if $F(t)$ is independent to $t$,
it results in confliction.
Therefore,
the solution \eqref{sol01} does not satisfies \eqref{eq2}.
{\hfill{$\square$}

For Eq. \eqref{sol02}, we have the following theorem.
\vspace*{-7pt}
\begin{theorem}\label{thm2}
The solution \eqref{sol02} does not satisfy system \eqref{eq2}.
\end{theorem}

{\bf Proof:}
It is easy to see that $x(t)$ satisfies   system \eqref{eq2} on $t\in [0,t_1]$.

For $t\in (t_1,T]$, by \eqref{sol02},
\begin{align*}
  x(t)&=E_{\frac{2}{3},1}(\rho t^{\frac{2}{3}})x_0
      +E_{\frac{2}{3},1}(\rho (t-t_1)^{\frac{2}{3}})y_1\\
      &+\int_{0}^t(t-\theta)^{-\frac{1}{3}}E_{\frac{2}{3},
      \frac{2}{3}}(\rho(t-\theta)^{\frac{2}{3}})\theta d\theta,
      \;\;t\in(t_1,T],
\end{align*}
we have\vspace*{-8pt}
\begin{align}\label{2}
 _0D_t^{\frac{2}{3}} x(t)
 &={}_0D_t^{\frac{2}{3}}E_{\frac{2}{3},1}(\rho t^{\frac{2}{3}})x_0
 +{}_0D_t^{\frac{2}{3}}E_{\frac{2}{3},1}(\rho (t-t_1)^{\frac{2}{3}})y_1\nonumber\\
 &+~_0D_t^{\frac{2}{3}}\int_{0}^t(t-\theta)^{-\frac{1}{3}}
 E_{\frac{2}{3},\frac{2}{3}}(\rho(t-\theta)^{\frac{2}{3}})\theta d\theta \nonumber\\
 &=\rho\Big[E_{\frac{2}{3},1}(\rho t^{\frac{2}{3}})x_0+\int_{0}^t(t-\theta)^{-\frac{1}{3}}E_{\frac{2}{3},\frac{2}{3}}
 (\rho(t-\theta)^{\frac{2}{3}})\theta d\theta\Big]+t\nonumber\\
 &+~_0D_t^{\frac{2}{3}} E_{\frac{2}{3},1}(\rho (t-t_1)^{\frac{2}{3}})y_1,
\end{align}
where\vspace*{-8pt}
\begin{align}
 &_0D_t^{\frac{2}{3}} E_{\frac{2}{3},1}(\rho (t-t_1)^{\frac{2}{3}})y_1
 =\int_0^t\frac{(t-s)^{-\frac{2}{3}}}{\Gamma(\frac{1}{3})}E'_{\frac{2}{3},1}(\rho (s-t_1)^{\frac{2}{3}})y_1ds\nonumber\\
 &=y_1\Big[\int_0^{t_1}\frac{(t-s)^{-\frac{2}{3}}}{\Gamma(\frac{1}{3})}E'_{\frac{2}{3},1}(\rho (s-t_1)^{\frac{2}{3}})ds
 +\int_{0}^{t-t_1}\frac{(t-t_1-s)^{-\frac{2}{3}}}{\Gamma(\frac{1}{3})}E'_{\frac{2}{3},1}(\rho s^{\frac{2}{3}})ds\Big]\nonumber\\
 &=y_1\Big[\int_0^{t_1}\frac{(t-s)^{-\frac{2}{3}}}{\Gamma(\frac{1}{3})}E'_{\frac{2}{3},1}(\rho (s-t_1)^{\frac{2}{3}})ds+
 ~_0D_{t-t_1}^{\frac{2}{3}} E_{\frac{2}{3},1}(\rho t^{\frac{2}{3}})\Big],\nonumber\\
 &=y_1\Big[\int_0^{t_1}\frac{(t-s)^{-\frac{2}{3}}}{\Gamma(\frac{1}{3})}E'_{\frac{2}{3},1}(\rho (s-t_1)^{\frac{2}{3}})ds+\rho E_{\frac{2}{3},1}(\rho (t-t_1)^{\frac{2}{3}})\Big].
\end{align}
Let $$y_1\int_0^{t_1}\frac{(t-s)^{-\frac{2}{3}}}{\Gamma(\frac{1}{3})}E'_{\frac{2}{3},1}(\rho (s-t_1)^{\frac{2}{3}})ds=G(t),$$
then
\begin{equation}\label{D0E}
_0D_t^{\frac{2}{3}} E_{\frac{2}{3},1}(\rho (t-t_1)^{\frac{2}{3}})y_1=\rho E_{\frac{2}{3},1}(\rho (t-t_1)^{\frac{2}{3}})y_1+G(t).
\end{equation}
Substituting \eqref{D0E} into \eqref{2} yields
\begin{align*}
 _0D_t^{\frac{2}{3}} x(t)
 &=\rho\Big[E_{\frac{2}{3},1}(\rho t^{\frac{2}{3}})x_0
      +E_{\frac{2}{3},1}(\rho(t-t_1)^{\frac{2}{3}})y_1\\
      &+\int_{0}^t(t-\theta)^{-\frac{1}{3}}
      E_{\frac{2}{3},\frac{2}{3}}(\rho(t-\theta)^{\frac{2}{3}})\theta d\theta\Big]+t+G(t)\\
      &=\rho x(t)+t+G(t).
\end{align*}
Similarly, if $G(t)\equiv 0$ for every $t\in (t_1,T]$, $G(t)$ is independent to $t$, which results in confliction.
Therefore,
\begin{align*}
 _0D_t^{\frac{2}{3}} x(t)
 &\neq \rho x(t)+t,~t\in(t_1,T].
\end{align*}
That is to say, Eq. \eqref{sol02} does not satisfy system \eqref{eq2}.
{\hfill{$\square$}

\begin{remark}\label{remA}
In view of the results of special case corresponding to Eq.\eqref{eq1} in Theorems \ref{thm1} and \ref{thm2}, for arbitrary continuous function $f(t)$,
$0<\alpha<1$, $i=1,2,\ldots,m,$
we can deduce that
\begin{equation}\label{1E-2}
_0D_t^\alpha\Big[\int_{t_i}^t(t-\theta)^{\alpha-1}E_{\alpha,\alpha}
(\rho(t-\theta)^\alpha)f(\theta)d\theta\Big]\neq \rho\Big[\int_{t_i}^{t}(t-\theta)^{\alpha-1}E_{\alpha,\alpha}
(\rho(t-\theta)^\alpha)f(\theta)d\theta\Big]+f(t),
\end{equation}
\begin{equation}\label{2E}
_0D_t^\alpha E_{\alpha,1}(\rho (t-t_i)^\alpha)y_i\neq \rho E_{\alpha,1}(\rho (t-t_i)^\alpha)y_i.
\end{equation}

If we replace $\rho$ by sectorial operator $A$, here $A$ is the infinitesimal generator of $\alpha$-resolvent operator $S_\alpha(t), T_\alpha(t)$ in Banach space $X$,
and take resolvent operator $S_\alpha(t)=E_{\alpha,1}(At^\alpha), T_\alpha(t)=t^{\alpha-1}E_{\alpha,\alpha}(At^\alpha)$ (see \cite{1}).
According to \eqref{1E-2} and \eqref{2E}, for $0<\alpha<1$, $i=1,2,\ldots,m,$ we obtain the following two results immediately,
\begin{equation}\label{D0T}
 _0D_t^\alpha\Big[\int_{t_i}^tT_\alpha(t-\theta)f(\theta)d\theta\Big]\neq A\Big[\int_{t_i}^{t}T_\alpha(t-\theta)f(\theta)d\theta\Big]+f(t),
\end{equation}
and
\begin{equation}\label{D0S}
_{0}D_t^\alpha S_\alpha(t-t_i)y_i\neq AS_\alpha(t-t_i)y_i.
\end{equation}
\end{remark}

Next, we prove the fact that $_{t_i}D_t^\alpha S_\alpha(t-t_i)y_i=AS_\alpha(t-t_i)y_i$.
\begin{lemma}\label{lem2}
Suppose $A$ is a sectorial operator of type $(M,\theta,\alpha,\mu)$.
If $0<\alpha<1, t>t_i,$ then
\[
  _{t_i}D_t^\alpha S_\alpha(t-t_i)y_i=AS_\alpha(t-t_i)y_i.
\]
\end{lemma}

{\bf Proof:}
From the definition of the Caputo derivative, we have
\begin{align*}
  _{t_i}D_{t}^{\alpha}S_\alpha(t-t_i)y_i
  &=\int_{t_i}^t \frac{(t-s)^{-\alpha}}{\Gamma(1-\alpha)}S'_{\alpha}(s-t_i)y_ids \\
  &=\int_{0}^{t-t_i} \frac{(t-t_i-s)^{-\alpha}}{\Gamma(1-\alpha)}S'_{\alpha}(s)y_ids
  ={}_{0}D_{t-t_i}^{\alpha}S_\alpha(t)y_i.
\end{align*}
By using the Laplace transform, combined with properties of the Laplace transform, we have
$$\mathcal{L}\{(t-t_i)^{-\alpha};\lambda\}=e^{-\lambda t_i}\Gamma(1-\alpha)\lambda^{-(1-\alpha)},~~
\mathcal{L}\{S_\alpha(t-t_i)y_i;\lambda\}=e^{-\lambda t_i}\lambda^{\alpha-1}R(\lambda^\alpha,A)y_i,$$
\vspace*{-16pt}
\begin{align}
\mathcal{L}\{{}_{0}D_{t-t_i}^{\alpha}S_\alpha(t)y_i;\lambda\}
&=\mathcal{L}\big\{\frac{(t-t_i)^{-\alpha}}{\Gamma(1-\alpha)}*S'_\alpha(t);\lambda\big\}y_i\nonumber\\
&=e^{-\lambda t_i}
\Big[\lambda^\alpha\mathcal{L}(S_\alpha(t);\lambda)-\lambda^{\alpha-1}\Big]y_i \nonumber\\
&=e^{-\lambda t_i}\lambda^{\alpha-1}(\lambda^\alpha I-A)^{-1}[\lambda^{\alpha}-(\lambda^{\alpha}-A)]y_i\nonumber\\
&=Ae^{-\lambda t_i}\lambda^{\alpha-1}R(\lambda^{\alpha},A)y_i
=A\mathcal{L}\{S_\alpha(t-t_i)y_i;\lambda\}.
\end{align}
Then we can accomplish the proof by taking the inverse Laplace transform.
{\hfill{$\square$}

\section{Main result \label {Sec4}}
In this section, we will deduce a new definition of mild solutions for problem \eqref{eq1} and justify the correctness of the new mild solution.
At first, using Remark \ref{remA}, Lemma \ref{lem2} obtained in Section \ref{Sec3}, we will prove that \eqref{sol1} and \eqref{sol2} are not solutions to system \eqref{eq1}.
Then by the method of undetermined coefficients, we will
get the conclusion of Theorem \ref{thm3}.
Finally, we will give a proof of Theorem \ref{thm3}.

\begin{theorem}\label{wrong}
The solution \eqref{sol1} and \eqref{sol2}  do not satisfy system \eqref{eq1}.
\end{theorem}

{\bf Proof:}
(i) By Lemma \ref{lem1}, for the expression of mild solution \eqref{sol1},  we can get that when $t\in [0,t_1],$
\begin{align*}
  _{0}D_{t}^{\alpha}x(t) &=~_{0}D_{t}^\alpha\Big[ S_\alpha(t)x_0+\int_{0}^tT_\alpha(t-\theta)f(\theta,x(\theta))d\theta\Big] \\
  & = A\Big[S_\alpha(t)x_0+\int_{0}^tT_\alpha(t-\theta)f(\theta,x(\theta))d\theta\Big]+f(t,x(t))\\
  &= Ax(t)+f(t,x(t)),
\end{align*}
where $x(t)$ satisfies system \eqref{eq1} on $[0,t_1]$.

For $t\in J_k$,
\[ x(t)= S_\alpha(t-t_k)[x(t_k^-)+I_k(x(t_k^-))]
          +\int_{t_k}^tT_\alpha(t-\theta)f(\theta,x(\theta))d\theta,\;t\in(t_k,t_{k+1}].\]
Applying the Caputo fractional derivative operator to $x(t)$ yields
\begin{align*}
 _0D^\alpha_t x(t)&=~_0D^\alpha_t\Big[S_\alpha(t-t_k)[x(t_k^-)+I_k(x(t_k^-))]\Big]+
 ~_0D^\alpha_t\Big[\int_{t_k}^tT_\alpha(t-\theta)f(\theta,x(\theta))d\theta\Big]\\
  &\ne AS_\alpha(t-t_k)[x(t_k^-)+I_k(x(t_k^-))]+A\int_{t_k}^tT_\alpha(t-\theta)f(\theta,x(\theta))d\theta+f(t,x(t)) ~~(\mbox{by Remark}~\ref{remA})\\
  &= Ax(t)+f(t,x(t)).
\end{align*}

(ii) Similarly, for the expression \eqref{sol2}, when $t\in J_k$, we have
\begin{align*}
_{0}D_{t}^{\alpha}x(t)&=~_{0}D_{t}^{\alpha} S_\alpha(t)x_0
+~_{0}D_{t}^{\alpha}\Big[\sum\limits_{i=1}^k S_\alpha(t-t_i)I_i(x(t_i^-))\Big]
+~_{0}D_{t}^{\alpha}\Big[\int_{0}^tT_\alpha(t-\theta)f(\theta,x(\theta))d\theta\Big]\\
&=AS_\alpha(t)x_0+\sum\limits_{i=1}^k~_{0}D_{t}^{\alpha} S_\alpha(t-t_i)I_i(x(t_i^-))
+A\int_{0}^tT_\alpha(t-\theta)f(\theta,x(\theta))d\theta+f(t,x(t))\\
&\neq AS_\alpha(t)x_0+\sum\limits_{i=1}^k AS_\alpha(t-t_i)I_i(x(t_i^-))
+A\int_{0}^tT_\alpha(t-\theta)f(\theta,x(\theta))d\theta+f(t,x(t))
~(\mbox{by \eqref{D0S}})\\
&=Ax(t)+f(t,x(t)).
\end{align*}

The above discission implies that \eqref{sol1} and \eqref{sol2} do not satisfy  system \eqref{eq1}.
{\hfill{$\square$}

\begin{remark}
From  the results of \eqref{DtT}, and Lemma \ref{lem2}, it is not difficult to obtain that if a function $x(t)\in PC(J,X)$ defined by the equation \eqref{sol1} is a mild solution for the equation
\begin{equation}\label{eqtk}
\left\{ \begin{array}{ll}
_{t_k}D_t^\alpha x(t)=Ax(t)+f(t,x(t)),\;\; 0<\alpha<1,\;t\in (t_k,t_{k+1}),\;k=0,1,\ldots,m,\\
\Delta x|_{t=t_k}=I_k(x(t_k^-)),\;k=1,2,\ldots,m,\\
x(0)=x_0.\\
\end{array} \right.
\end{equation}
\end{remark}

In what follows, we deduce the formula of mild solution for system \eqref{eq1} through the method of undetermined coefficients.

(i) When $t\in [0,t_1]$, the mild solution is (see \cite{2})
$$x(t)=S_\alpha(t)x_0+\int_0^tT_\alpha(t-s)f(s)ds,\ t\in [0,t_1].$$

(ii) For $t\in (t_1,t_2]$, suppose that
$$x(t)=S_\alpha(t)z+\int_0^tT_\alpha(t-s)f(s)ds,\ \ z\in X.$$
By the impulsive condition
$$x(t_1^+)-x(t_1^-)=S_\alpha(t_1)[z-x_0]=I_1(x_1^-),$$
we obtain that
$$z=x_0+S^{-1}_\alpha(t_1)I_1(x_1^-). $$

Thus,
$$x(t)=S_\alpha(t)x_0+S_\alpha(t)S^{-1}_\alpha(t_1)I_1(x_1^-)+
\int_0^{t}T_\alpha(t-s)f(s)ds,\ t\in (t_1,t_2].$$

By the same way, for every $t\in (t_k, t_{k+1}]$, we have
$$x(t)=S_\alpha(t)x_0+\sum_{i=1}^kS_\alpha(t)S^{-1}_\alpha(t_i)I_i(x_i^-)+
\int_0^{t}T_\alpha(t-s)f(s)ds,\  t\in (t_k, t_{k+1}].$$

\begin{remark}\label{inverse}
When $0<\alpha<1$, $S_\alpha(t)$ is bounded linear and positive operator in Banach space $X$.
It is easy to see $S_\alpha(t)$ is a one-to-one mapping. So, for each impulsive time $t_i,(i=1,2,\ldots m)$, the inverse operator $S^{-1}_\alpha(t_i)$ exists.
\end{remark}


{\bf Proof of Theorem 1.1:}

When $t\in [0,t_1]$,
the mild solution is
$$x(t)=S_\alpha(t)x_0+\int_0^tT_\alpha(t-s)f(s)ds,\ t\in [0,t_1].$$

For $t\in (t_k,t_{k+1}]$, by Lemma \ref{lem1},
we get
\begin{align*}
 _0D^\alpha_t x(t)& =~_0D^\alpha_t\Big(S_\alpha(t)x_0+S_\alpha(t)\sum\limits_{i=1}^k S^{-1}_\alpha(t_i)I_i(x(t_i^-))+
 \int_0^tT_\alpha(t-\theta)f(\theta,x(\theta))d\theta\Big)\\
  &=AS_\alpha(t)x_0+A S_\alpha(t)\sum_{i=1}^k S^{-1}_\alpha(t_i)I_i(x(t_i^-))+A\int_0^tT_\alpha(t-\theta)f(\theta,x(\theta))d\theta+f(t,x(t))\\
  &=A\Big[ S_\alpha(t)x_0+S_\alpha(t)\sum_{i=1}^k S^{-1}_\alpha(t_i)I_i(x(t_i^-))+\int_0^tT_\alpha(t-\theta)f(\theta,x(\theta))d\theta\Big]+f(t,x(t))\\
  &=Ax(t)+f(t,x(t)).
\end{align*}

For $t=0$, we have $$x(0)=S_\alpha(0)x_0+\int_0^0T_\alpha(0-\theta)f(\theta,x(\theta))d\theta=x_0.$$

Moreover, for $t=t_k$, we have
\begin{align*}
  \Delta x(t_k) &=x(t_k^+)-x(t_k^-)
   =S_\alpha(t_k)\sum\limits_{i=1}^k S^{-1}_\alpha(t_i)I_i(x(t_i^-))
   -S_\alpha(t_k)\sum\limits_{i=1}^{k-1}S^{-1}_\alpha(t_i)I_i(x(t_i^-))\\
  & =S_\alpha(t_k)S^{-1}_\alpha(t_k)I_k(x(t_k^-))
   =I_k(x(t_k^-)).
\end{align*}
It thus follows  that  \eqref{sol3} is a solution to the linear impulsive fractional differential Eq.\eqref{eq1}.
The proof of Theorem 1.1 is complete.
{\hfill{$\square$}

\end{document}